\documentclass[12pt]{article}
\usepackage{amssymb,amsmath,amsthm,hyperref}
\usepackage{graphicx}
\graphicspath{ {./} }

\usepackage[utf8]{inputenc}
\usepackage{amsthm}
\usepackage{amsmath}
\usepackage{amssymb}
\usepackage{enumitem}
\setlist{topsep = 0pt}
\usepackage[table]{xcolor}
\usepackage{tikz}
\usepackage{xcolor}
\usepackage{changepage}

\theoremstyle{plain}
\newtheorem{thm}{Theorem}[section]

\theoremstyle{definition}

\newtheorem{rmk}[thm]{Remark}
\newtheorem{exmp}[thm]{Example}
\usepackage{caption}
\usepackage{subcaption}
\usepackage{tikz}
\usetikzlibrary{arrows}

\usepackage{graphicx}
\usepackage{xcolor} 

\usepackage{dirtytalk} 
\usepackage{mathrsfs} 

\usepackage{esint} 
\usepackage{marvosym} 

\usepackage{titlefoot}
\usepackage{titlesec}
\titleformat*{\section}{\centering\MakeUppercase}
\titleformat*{\subsection}{\normalsize\bfseries}
\titlespacing\section{0pt}{6pt plus 4pt minus 2pt}{0pt plus 2pt minus 2pt}
\usepackage[margin=2cm]{geometry}
\newcommand{\comm}[1]{}

\usepackage[capitalise]{cleveref}
\definecolor{pu}{rgb}{0.56, 0, 0.56}
\definecolor{gr}{rgb}{0, 0.67, 0}

\begin{document}
\title{Implementation: The conjugacy problem in right-angled Artin groups}
\author{\small{\uppercase{Gemma Crowe, Michael Jones}}}
\date{}
\maketitle
\begin{abstract}
   In 2009, Crisp, Godelle and Wiest constructed a linear-time algorithm to solve the conjugacy problem in right-angled Artin groups. This algorithm has now been implemented in Python, and the code is freely available on \href{https://github.com/gmc369/Conjugacy-problem-RAAGs}{GitHub}. This document provides a summary of how the code works. As well as determining whether two elements $w_{1}, w_{2}$ are conjugate in a RAAG $A_{\Gamma}$, our code also returns a conjugating element $x \in A_{\Gamma}$ such that $w_{1} = x^{-1}w_{2}x$, if $w_{1}$ and $w_{2}$ are conjugate.
\end{abstract}
\section{Introduction}
In their paper `\textit{The conjugacy problem in subgroups of right-angled Artin groups}', Crisp, Godelle and Wiest created a linear-time algorithm to solve the conjugacy problem in right-angled Artin groups (RAAGs) \cite{CGW}. This algorithm has now been implemented as a Python program, and is freely available on \href{https://github.com/gmc369/Conjugacy-problem-RAAGs}{GitHub}.
\par 
This document provides an overview of the code. We recommend the reader first takes a look at the original paper \cite{CGW}, to understand the key tools and steps of this algorithm. 
\begin{rmk}
Our Python code requires the \texttt{networkx} module. You may also need to install Pillow and nose. Details on how to install these modules can be found here: \url{https://networkx.org/documentation/stable/install.html}.
\end{rmk}

\subsection{Notation}
Throughout we let $A_{\Gamma}$ denote a RAAG with defining graph $\Gamma$. We let $N$ denote the number of vertices in the generating set, i.e. the size of the standard generating set for $A_{\Gamma}$. We assume generators commute if and only if there does not exist an edge between corresponding vertices in $\Gamma$, to match with the convention used in \cite{CGW}. 
\par
Examples are provided both in the code as well as this document to assist the user. Throughout this summary, we will often use the RAAG defined in Example 2.4 from \cite{CGW}, which has the following presentation:
\begin{equation}\label{CGW example}
    A_{\Gamma} = \langle a_{1}, a_{2}, a_{3}, a_{4} \; | \; [a_{1}, a_{4}] = 1, [a_{2}, a_{3}] = 1, [a_{2}, a_{4}] = 1 \rangle.
\end{equation}
This RAAG is defined by the following graph:
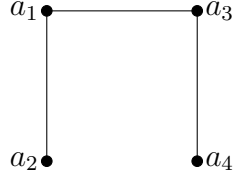
\begin{figure}
    \centering
\begin{tikzpicture}
    \filldraw[black] (0,0) circle (2pt);
    \filldraw[black] (2,0) circle (2pt);
    \filldraw[black] (2,2) circle (2pt);
    \filldraw[black] (0,2) circle (2pt);
    \draw (0,0) -- (0,2);
    \draw (0,2) -- (2,2);
    \draw (2,2)-- (2,0);
    \node at (-0.3, 0) {$a_{2}$};
    \node at (-0.3, 2) {$a_{1}$};
    \node at (2.3, 0) {$a_{4}$};
    \node at (2.3, 2) {$a_{3}$};
    \end{tikzpicture}
    \caption{Defining graph $\Gamma$}
    \label{fig:my_label}
\end{figure}

\section{Summary of algorithm}
We summarise the key steps of the algorithm from \cite{CGW}. For further details on piling representations, see Chapter 14 of \cite{OfficeHour}.
\par 
\underline{ALGORITHM}: Conjugacy problem in RAAGs. 
\par 
\textbf{Input}:
\begin{enumerate}
    \item RAAG $A_{\Gamma}$: number of vertices of the defining graph, and list of commuting generators.
    \item Words $v,w$ representing group elements in $A_{\Gamma}$.
\end{enumerate}
\textbf{Step 1: Cyclic reduction}
\begin{adjustwidth}{1.5cm}{}
    Produce the piling representation $\pi(v)$ of the word $v$, and apply cyclic reduction to $\pi(v)$ to produce a cyclically reduced piling $p$. Repeat this step for the word $w$ to get a cyclically reduced piling $q$.
\end{adjustwidth}
\textbf{Step 2: Factorisation}
\begin{adjustwidth}{1.5cm}{}
    Factorise each of the pilings $p$ and $q$ into non-split factors. If the collection of subgraphs do not coincide, \textbf{Output} = \texttt{\color{blue}False}.\\
    Otherwise, continue to Step 3.
\end{adjustwidth}
\textbf{Step 3: Compare non-split factors}
\begin{adjustwidth}{1.5cm}{}
    If $p = p^{(1)}\dots p^{(k)}$ and $q = q^{(1)}\dots q^{(k)}$ are the factorisations found in Step 2, then for each $i = 1, \dots, k$, do the following:
    \begin{adjustwidth}{1.5cm}{}
        \begin{enumerate}[label=(\roman*)]
            \item Transform the non-split cyclically reduced pilings $p^{(i)}$ and $q^{(i)}$ into pyramidal pilings $\Tilde{p}^{(i)}$ and $\Tilde{q}^{(i)}$.
            \item Produce the words representing these pilings in cyclic normal form $\sigma^{\ast}\left(\Tilde{p}^{(i)}\right)$ and $\sigma^{\ast}\left(\Tilde{q}^{(i)}\right)$.
            \item Decide whether $\sigma^{\ast}\left(\Tilde{p}^{(i)}\right)$ and $\sigma^{\ast}\left(\Tilde{q}^{(i)}\right)$ are equal up to a cyclic permutation. If not, \textbf{Output} = \texttt{\color{blue}False}.
        \end{enumerate}
    \end{adjustwidth}
    \textbf{Output} = \texttt{\color{blue}True}.
\end{adjustwidth}
    
\begin{exmp}
We provide an example of how to implement this algorithm with our Python code, using the RAAG $A_{\Gamma}$ defined in Equation \ref{CGW example}. Consider the following two words:
\[ w_{1} = a^{-2}_{2}a^{-1}_{4}a_{3}a_{2}a_{4}a_{1}a_{2}a^{-1}_{1}a^{2}_{2}a^{-1}_{4}, \quad w_{2} = a_{4}a_{3}a^{-1}_{4}a_{2}a_{1}a_{2}a^{-1}_{1}a^{-1}_{4}.
\]
These represent conjugate elements in $A_{\Gamma}$, since
\[ w_{1} = \left(a^{2}_{4}a^{2}_{2}\right)^{-1}\cdot w_{2} \cdot  \left(a^{2}_{4}a^{2}_{2}\right).
\]
To check this, we input the following code. See \cref{sec:inputs} on how to input words and commutation relations from a RAAG.
\begin{align*}
    \texttt{w\_{1}} &= \texttt{[-2,-2,-4,3,2,4,1,2,-1,2,2,-4]}\\
    \texttt{w\_{2}} &= \texttt{[4,3,-4,2,1,2,-1,-4]} \\
    \texttt{CGW\_Edges} &= \texttt{[(1,4), (2,3), (2,4)]}\\
    &\texttt{is\_conjugate(w\_1, w\_2, 4, CGW\_Edges)}
\end{align*}
The output from this is
\[ \texttt{\color{blue}True}, \texttt{[-2, -2, -4, -4]}
\]
The first argument is a \texttt{\color{blue}True} or \texttt{\color{blue}False} statement of whether $w_{1}$ and $w_{2}$ are conjugate. The second argument returns a conjugating element $x$ such that $w_{1} = x^{-1}w_{2}x$. 
\end{exmp}
It is important to note that the order of the input words determines the conjugating element. In particular, when we input $w_{1}$ following by $w_{2}$, we obtain a conjugator $x$ such that $w_{1} = x^{-1}w_{2}x$. If we reverse the order of $w_{1}$ and $w_{2}$, the program will return the inverse of this conjugator, since $w_{2} = xw_{1}x^{-1}$. 

\section{Inputs}\label{sec:inputs}
\subsection{Words}
The approach of using symbols $a,b,c...$ to represent letters in a word is limiting, since there are only finitely many letters in the alphabet. Instead, we use the positive integers to represent the generators, and the negative integers to represent their inverses:
\begin{center}
\begin{tabular}{|c|c||c|c|}
	\hline
	Generator & Representation & Generator & Representation\\
	\hline\hline
	$a$ & \texttt{1} & $a^{-1}$ & \texttt{-1}\\
	\hline
	$b$ & \texttt{2} & $b^{-1}$ & \texttt{-2}\\
	\hline
	$c$ & \texttt{3} & $c^{-1}$ & \texttt{-3}\\
	\hline
	... & ... & ... & ... \\
	\hline
\end{tabular}
\end{center}
This also matches with our convention of ordering which will be shortlex, i.e.
\[ 1 < -1 < 2 < -2 < \dots 
\]
We note this convention is the opposite of the normal form convention in \cite{CGW}. 
\par 
Words in $A_{\Gamma}$ are represented by lists of integers. The following table gives some examples. Here $\varepsilon$ denotes the empty word.
\begin{center}
    \begin{tabular}{|c|c|}
	\hline
	Word & Representation\\
	\hline\hline
	$abc$ & \texttt{[1,2,3]}\\
	\hline
	$\varepsilon$ & \texttt{[]}\\
	\hline
	$a^{-2}c^3b^{-1}$ & \texttt{[-1,-1,3,3,3,-2]}\\
	\hline
	$a^{100}b^{-1}$ & \texttt{([1]*100)+[-2]}\\
	\hline
	$abc...xyz$ & \texttt{\color{pu}list\color{black}(\color{pu}range\color{black}(1,27))}\\
	\hline
\end{tabular}
\end{center}

\subsection{Commuting generators}
In many of the functions in our code, we need to define which generators commute in $A_{\Gamma}$. This is represented by a list of tuples, each of which represents a pair of commuting generators. Here are some examples:
\begin{center}
\begin{tabular}{|c|c|}
	\hline
	Group & Commuting generators list\\
	\hline\hline
	$\mathbb{F}_n$ & \texttt{[]}\\
	\hline
	$\mathbb{Z}^2$ & \texttt{[(1,2)]}\\
	\hline
	$\mathbb{F}_2\times\mathbb{Z}$ & \texttt{[(1,3),(2,3)]}\\
	\hline
	$\mathbb{Z}^n$ & \texttt{[(i,j) \color{orange}for \color{black}i \color{orange}in \color{pu}range\color{black}(1,n) \color{orange}for \color{black}j \color{orange}in \color{pu}range\color{black}(i+1,n+1)]}\\
	\hline
\end{tabular}
\end{center}
By convention, we do not include edges $(n,n)$, and we order our list by shortlex ordering. For example, if $(3,4)$ is in our list, we do not need to also include $(4,3)$. Also, the code does not require inverses of generators to be taken into account. In particular, for every tuple $(n,m)$ in the list, we do not need to also include $(-n, m), (n, -m)$ or $(-n, -m)$.
\par 
This list should contain precisely one tuple for each non-edge in the defining graph. If no list is entered, the program takes \texttt{[]} as default (i.e. the free group).

\begin{rmk}
In the Python code, we use the command \texttt{commuting\_elements} to denote commuting generators.
\end{rmk}

\section{Pilings}
Pilings are represented by lists of lists. Each list within the main list corresponds to a column of the piling, reading from bottom to top. A `$+$' bead is represented by \texttt{1}, a `$-$' bead is represented by \texttt{-1} and a `$0$' bead is represented by \texttt{0}. Here are some examples:
\begin{center}
\begin{tabular}{|c|c|}
	\hline
	Piling & Representation\\
	\hline\hline
	\includegraphics[scale=0.1]{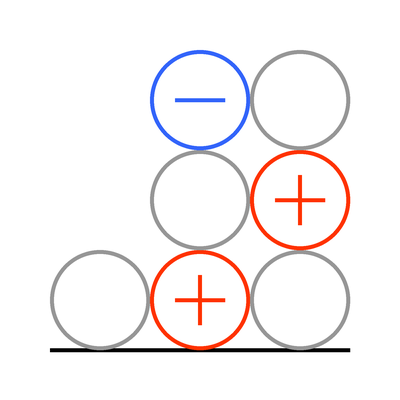} & \texttt{[[0],[1,0,-1],[0,1,0]]}\\
	\hline
	\includegraphics[scale=0.1]{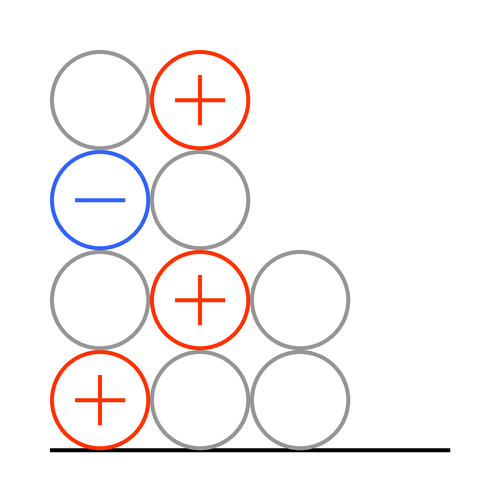} & \texttt{[[1,0,-1,0],[0,1,0,1],[0,0],[]]}\\
	\hline
	\includegraphics[scale=0.1]{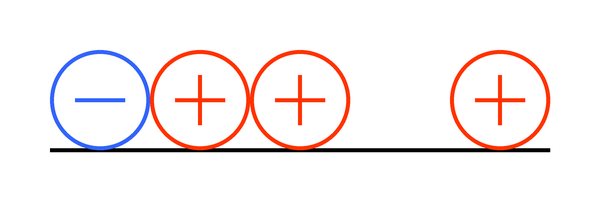} & \texttt{[[-1],[1],[1],[],[1]]}\\
	\hline
\end{tabular}
\end{center}
We note the empty piling is represented by \texttt{[[]$\ast N$]} where $N$ denotes the number of generators.
\par 
Reading pilings as a list of lists is not as user friendly as a pictorial representation. We refer the reader to Appendix \ref{appen:draw}, which explains how to use the \texttt{draw\_piling()} function. This allows the user to create a pictorial representation of pilings.

\subsection{Programs}
One of the most useful steps of the algorithm is converting words, which represent group elements in a RAAG, into piling representations. One key fact from these constructions is that if two words $u,v \in A_{\Gamma}$ represent the same group element, then the piling representations for $u$ and $v$ will be the same. In particular, every group element in $A_{\Gamma}$ is uniquely determined by its piling. We describe this function here.
\par 
\textbf{Function:} \texttt{piling(w, N, commuting\_elements=[])}
\par
\textbf{Input}:
\begin{enumerate}
    \item Word $w$ representing a group element in $A_{\Gamma}$.
    \item $N=$ number of vertices in defining graph.
    \item List of commuting generators.
\end{enumerate}
\textbf{Output}: Piling representation of $w$, as a list of lists.
\par 
By construction, the piling will reduce any trivial cancellations in $w$ after shuffling, so we do not require our input word $w$ to be reduced.
\begin{exmp}
Suppose you wish to compute the piling for the word $ab^{2}a^{-1}b\in\mathbb{F}_2$. Then you would input \texttt{piling([1,2,2,-1,2], 2)}, which outputs the following:
\[ \texttt{[[1, 0, 0, -1, 0], [0, 1, 1, 0, 1]]}
\]
Similarly we could compute the piling for the word $a^{-2}_{2}a^{-1}_{4}a_{3}a_{2}a_{4}a_{1}a_{2}a^{-1}_{1}a^{2}_{2}a^{-1}_{4} \in A_{\Gamma}$ from Equation \ref{CGW example}. 
This outputs the following piling:
\[ \texttt{[[0, 0, 1, 0, -1, 0, 0], [-1, 0, 1, 0, 1, 1], [0, 1, 0, 0], [-1, 0]]}
\]
\begin{figure}[ht]
    \centering
    \subfloat[\centering]{{\includegraphics[scale = 0.2]{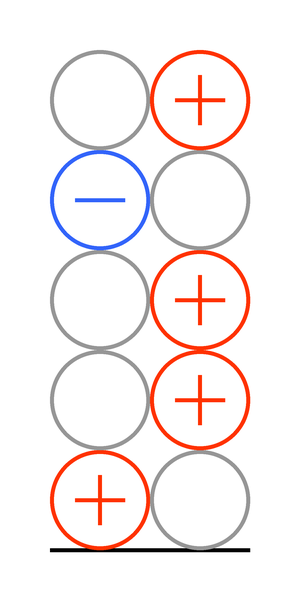} }}%
    \qquad
    \subfloat[\centering]{{\includegraphics[scale = 0.2]{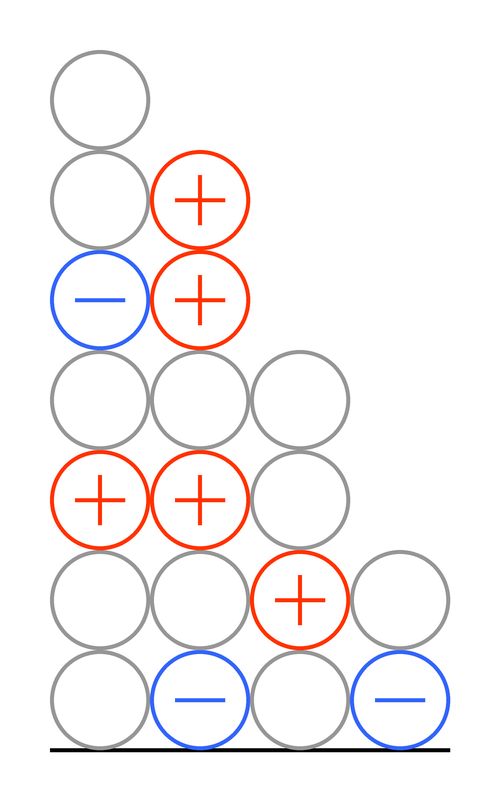} }}%
    \caption{$ab^{2}a^{-1}b\in\mathbb{F}_2$ and $a^{-2}_{2}a^{-1}_{4}a_{3}a_{2}a_{4}a_{1}a_{2}a^{-1}_{1}a^{2}_{2}a^{-1}_{4} \in A_{\Gamma}$ }%
    \label{fig:ex pilings}%
\end{figure}
Figure \ref{fig:ex pilings} gives a pictorial representation of these pilings. Each list represents a column in the piling.
\end{exmp}
We now present a function which computes the normal form representative from a piling. For us, this is the shortlex shortest word which represents the piling. We remind the reader that this convention is the opposite ordering for normal forms in \cite{CGW}.
\newpage
\textbf{Function:} {\texttt{word(p, commuting\_elements=[])}}
\par 
\textbf{Input}:
\begin{enumerate}
    \item $p$ = piling.
    \item List of commuting generators.
\end{enumerate}
\textbf{Output}: normal form word which represents the piling $p$.

\begin{exmp}
The following is an example of how to use this function for the following piling in $A_{\Gamma} = \langle a,b,c \; | \; ac=ca\rangle$:
\begin{align*}
	&\texttt{p=[[1,0],[0,0,-1],[-1,0]]}\\
	&\texttt{w=word(p,[(1,3)])}\\
	&\texttt{\color{pu}print\color{black}(w)}
\end{align*}
The output is:
\begin{center}
	\texttt{[1,-3,-2]}
\end{center}
\end{exmp}
\section{Cyclic reduction}
After constructing pilings from our input words, the next step of the conjugacy algorithm is to cyclically reduce each piling.
\par 
\textbf{Function:} {\texttt{cyclically\_reduce(p, commuting\_elements=[])}}
\par 
\textbf{Input}:
\begin{enumerate}
    \item $p =$ piling.
    \item List of commuting generators.
\end{enumerate}
\textbf{Output}:
\begin{enumerate}
    \item Cyclically reduced piling.
    \item Conjugating element.
\end{enumerate}
\begin{exmp}
Let $w = abca^{-1}\in\langle a,b,c\; | \; bc=cb\rangle$. We can compute all possible cyclic reductions on $w$ as follows:
\begin{align*}
	&\texttt{p=piling([1,2,3,-1], 3,[(2,3)])}\\
	&\texttt{p\_cr=cyclically\_reduce(p,[(2,3)])}\\
	&\texttt{w=word(p\_cr[0],[(2,3)])}\\
	&\texttt{\color{pu}print\color{black}(w, p\_cr[1])}
\end{align*}
The output from this is:
\begin{center}
	\texttt{[2,3], [1]}
\end{center}
This function is the first example where we output a conjugating element. In this example, when we cyclically reduce $w = a\cdot bc \cdot a^{-1}$ to the word $bc$, we have taken out the conjugating element $a$. 
\end{exmp}

\section{Graphs}
We remind the reader that for the functions in this section, we need to import the \texttt{networkx} module into Python. This makes working with the defining graph much easier, in particular computing induced subgraphs and connected components. 
\par
The next step in the conjugacy algorithm is to compute the non-split factors from a word based on the defining graph. First, we need a method to input our defining graph.
\par 
\textbf{Function:} {\texttt{graph\_from\_edges(N, commuting\_edges=[])}}
\par 
\textbf{Input}:
\begin{enumerate}
    \item $N =$ number of vertices in defining graph.
    \item List of commuting generators.
\end{enumerate}
\textbf{Output}: \texttt{networkx} graph, which  represents the defining graph.
\par 
We recall that in \cite{CGW}, the convention is to add edges for non-commuting vertices. This is precisely what this function does - the output is a graph on $N$ vertices, with edges between vertices $(x,y)$ if and only if $(x,y)$ is not in the list of commuting generators. 
\begin{exmp}\label{exmp: graph}
Let's compute the defining graph from Equation \ref{CGW example}. We input the following:
\begin{align*}
    &\texttt{CGW\_Edges = [(1,4), (2,3), (2,4)]} \\
    &\texttt{g = graph\_from\_edges(4, CGW\_Edges)} \\
    &\texttt{nx.draw(g)}
\end{align*}
This outputs the following graph:
\begin{figure}
    \centering
    \includegraphics[scale = 0.6]{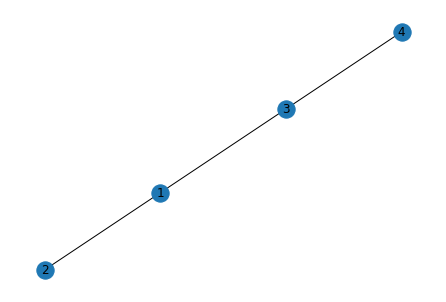}
    \caption{\texttt{networkx} graph defined in Example 6.1.}
    \label{fig:my_label2}
\end{figure}
\end{exmp}
\subsection{Factorising}
We now want to build up a function which takes our defining graph $g$ and our input word $w$, and computes the non-split factors related to $w$. 
\par 
\textbf{Function:} {\texttt{factorise(g, p)}}
\par
\textbf{Input}:
\begin{enumerate}
    \item $g =$ \texttt{networkx} graph representing the defining graph of $A_{\Gamma}$.
    \item $p =$ piling.
\end{enumerate}
\textbf{Output}: List of subgraphs of $g$, one for each non-split factor of $p$.
\par 
This function first checks which columns in the piling contain a \texttt{1} or \texttt{-1} bead (this is the \texttt{supp\_piling()} function). From this subset of vertices, we build the induced subgraph, and then return a list of the connected components of this subgraph. Each of these connected components represents the non-split factors.
\par 
At this stage, we have extracted each of the subgraphs which correspond to the factors. We now want to extract pilings which represent each factor. 
\comm{
\par 
\textbf{Function:} {\texttt{graphs\_to\_nsfactor(g, w)}}
\par 
\textbf{Input}:
\begin{enumerate}
    \item $g =$ connected component of defining graph.
    \item $w =$ word representing a group element of $A_{\Gamma}$.
\end{enumerate}
\textbf{Output}: non-split factor of $w$ based on $g$.}
\par 
\textbf{Function:} {\texttt{graphs\_to\_nsfactors(l, w, N, commuting\_elements)}}
\par 
\textbf{Input}:
\begin{enumerate}
    \item $l =$ list of subgraphs representing connected components.
    \item $w =$ word representing a group element of $A_{\Gamma}$.
    \item $N =$ number of vertices in defining graph.
    \item List of commuting generators.
\end{enumerate}
\textbf{Output}: list of the non-split factors from $w$ represented by pilings.
\begin{exmp}
We use the example from Figure 4 of \cite{CGW}. We input the following:
\begin{align*}
    &\texttt{w = [2,3,-4]}\\
    &\texttt{CGW\_Edges = [(1,4), (2,3), (2,4)]} \\
    &\texttt{g = graph\_from\_edges(4, CGW\_Edges)}\\
    &\texttt{graphs = factorise(g, piling(w, 4, CGW\_Edges))} \\
    &\texttt{\color{pu}print\color{black}(graphs\_to\_nsfactors(graphs, w, 4, CGW\_Edges))}
\end{align*}
The output is 
\[ \texttt{[[[0], [1], [], []], [[0], [], [1, 0], [0, -1]]]}
\]
where the first piling represents the factor $a_{2}$, and the second piling represents the factor $a_{3}a^{-1}_{4}$.
\end{exmp}

\section{Pyramidal Pilings}
For the final step in the algorithm, we need a function which converts pilings into pyramidal pilings. After this, we need to check when two pyramidal pilings are equal up to cyclic permutation.
\par 
\textbf{Function:} {\texttt{pyramidal\_decomp(p, commuting\_elements=[])}}
\par 
\textbf{Input}:
\begin{enumerate}
    \item $p = $ non-split piling.
    \item List of commuting generators.
\end{enumerate}
\textbf{Output}: factors $p_{0}, p_{1}$ as pilings.
\par 
It is important to note that the input piling must be non-split, otherwise we cannot achieve a pyramidal piling and the function will run forever. 
\begin{exmp}
Let's compute the decomposition in Figure 5 of \cite{CGW}. We input the following:
\begin{align*}
    &\texttt{CGW\_Edges = [(1,4), (2,3), (2,4)]} \\
    &\texttt{p = [[0,1,0,-1,0], [0,1,0,1], [0,1,0,0], [-1,0]]} \\
    &\texttt{decomp = pyramidal\_decomp(p, CGW\_Edges)}\\
    &\texttt{\color{pu}print\color{black}(decomp)}
\end{align*}
The output is then
\[ \texttt{([[0], [], [0, 1], [-1, 0]], [[1, 0, -1, 0], [0, 1, 0, 1], [0, 0], []])}
\]
\end{exmp}

\textbf{Function:} {\texttt{pyramidal(p, N, commuting\_elements=[])}}
\par 
\textbf{Input}:
\begin{enumerate}
    \item $p = $ non-split cyclically reduced piling.
    \item $N = $ number of vertices in defining graph.
    \item List of commuting generators.
\end{enumerate}
\textbf{Output}:
\begin{enumerate}
    \item Pyramidal piling of $p$.
    \item Conjugating element.
\end{enumerate}
The \texttt{pyramidal()} function iteratively applies \texttt{pyramidal\_decomp()} until the factor $p_{0}$ is empty. Since the only operation we apply on the piling is cyclic permutations, we can add these steps to our conjugating element.
\begin{exmp}
Again we take the same example above. Our input is:
\begin{align*}
    &\texttt{CGW\_Edges = [(1,4), (2,3), (2,4)]} \\
    &\texttt{p = [[0,1,0,-1,0], [0,1,0,1], [0,1,0,0], [-1,0]]} \\
    &\texttt{pyr = pyramidal(p, 4, CGW\_Edges)}
\end{align*}
Figure \ref{fig:pyr} gives the pictorial representation of this pyramidal piling. The conjugating element output in this example is \texttt{\color{blue}[-4,3,-4]}.
\begin{figure}[ht]
    \centering
    \includegraphics[scale = 0.2]{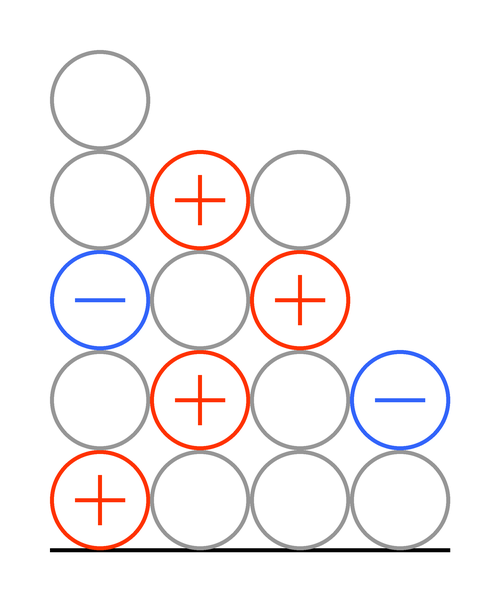}
    \caption{Pyramidal piling.}
    \label{fig:pyr}
\end{figure}
\end{exmp}
\vspace{-15pt}
\textbf{Function:} {\texttt{is\_cyclic\_permutation(w,v)}}
\par 
\textbf{Input}: Two words $w,v$ representing group elements in $A_{\Gamma}$.\\
\textbf{Output}:
If $w$ is a cyclic permutation of $v$, the function outputs two arguments:
\begin{enumerate}
    \item \texttt{\color{blue}True}
    \item Conjugating element.
\end{enumerate}
Otherwise, return \texttt{\color{blue}False}.

\section{Implementation of the Conjugacy Problem}
\subsection{Simple Cases}
If we want to solve the conjugacy problem in either $\mathbb{F}_n$ or $\mathbb{Z}^n$, the algorithm for implementing the conjugacy problem is more straightforward, and using the \texttt{is\_conjugate()} function detailed below is wasteful. Certainly in the latter case, the conjugacy problem is trivial; it is equivalent to checking if two words are equal, which is just a case of comparing the number of occurrences of each letter in the words. For free groups, two cyclically reduced words are conjugate if and only if they are cyclic permutations of each other. Hence one can make a shorter program for free groups by simply using the \texttt{cyclically\_reduce()} and \texttt{is\_cyclic\_permutation()} functions.
\par
For general RAAGs, we now have the necessary functions to implement the conjugacy problem.
\par 
\textbf{Function:} {\texttt{is\_conjugate(w1, w2, N, commuting\_elements=[])}}
\par 
\textbf{Input}:
\begin{enumerate}
    \item Two words $w_{1}, w_{2}$ representing group elements in $A_{\Gamma}$.
    \item $N = $ number of vertices in defining graph.
    \item List of commuting generators.
\end{enumerate}
\textbf{Output}: If $w_{1}$ is conjugate to $w_{2}$ in $A_{\Gamma}$ then return:
\begin{enumerate}
    \item \texttt{\color{blue}True}
    \item Conjugating element $x$ such that $w_{1} = x^{-1}w_{2}x$.
\end{enumerate}
Otherwise, return \texttt{\color{blue}False}.
\par 
We note that when computing the conjugating element, we find a reduced word $x$ which represents an element which conjugates $w_{1}$ to $w_{2}$.

\appendix
\section{The `draw piling' Function}\label{appen:draw}
The representation of pilings as a list of lists of \texttt{1}s, \texttt{-1}s and \texttt{0}s is not very readable in the Python interpreter. For small pilings it is tolerable, but for larger ones the \texttt{draw\_piling()} function is useful.
\par 
\noindent\texttt{draw\_piling()} can take a total of nine arguments, however only the first one, which is the piling to be drawn, is needed. The function returns nothing, but by default it will show a picture of the piling in a new window, as well as save a \texttt{PNG} file of it. 
\par 
The following table lists all the arguments of the \texttt{draw\_piling()} function and what they do:
\begin{center}
\begin{tabular}{|c|c|c|}
	\hline
	Argument & Type & Description\\
	\hline\hline
	\texttt{piling} & Piling & \tiny The (mandatory) piling to draw.\\
	\hline
	\texttt{scale} & Float & \tiny The scale to draw the piling at. Default is \texttt{100.0}.\\
	\hline
	\texttt{plus\_colour} & Colour & \tiny The colour to draw the `$+$' beads. Default is red.\\
	\hline
	\texttt{zero\_colour} & Colour & \tiny The colour to draw the `$0$' beads. Default is grey.\\
	\hline
	\texttt{minus\_colour} & Colour & \tiny The colour to draw the `$-$' beads. Default is blue.\\
	\hline
	\texttt{anti\_aliasing} & Integer & \tiny The super-sampling resolution for anti-aliasing.\\
	& & \tiny Only allows positive integers. Default is 4.\\
	\hline
	\texttt{filename} & String & \tiny The filename to save the piling as. Default is \texttt{\color{gr}"piling.png"}.\\
	\hline
	\texttt{show} & Boolean & \tiny Whether to show the piling in a window. Default is \texttt{\color{orange}True}.\\
	\hline
	\texttt{save} & Boolean & \tiny Whether to save the piling. Default is \texttt{\color{orange}True}.\\
	\hline
\end{tabular}
\end{center}

\section{Solution to the Word Problem}
With the tools that \texttt{pilings.py} provides, it is not hard to solve the Word Problem in any given RAAG. The following code defines a function that decides if a given word is equal to the identity:
\begin{align*}
	\texttt{\color{orange}def} \; &\texttt{\color{blue}identity\color{black}(w, N, commuting\_elements=[]):}\\
	&\texttt{\indent p=piling(w, N, commuting\_elements)\color{red}\#generate reduced piling}\\
	&\texttt{\indent reduced\_w=word(p, commuting\_elements)\color{red}\#read off reduced word}\\
	&\texttt{\indent \color{orange}return\color{black}(reduced\_w==[])\color{red}\#return whether it equals the identity}
\end{align*}
The following code defines a function that decides if two given words are equal in a given RAAG:
\begin{align*}
	\texttt{\color{orange}def} \; &\texttt{\color{blue}equal\color{black}(w1, w2, N, commuting\_elements=[]):}\\
	&\texttt{\indent \color{red}\#generate reduced pilings}\\
	&\texttt{\indent p1=piling(w1, N, commuting\_elements)}\\
	&\texttt{\indent p2=piling(w2, N, commuting\_elements)}\\
	&\texttt{\indent \color{red}\#read off reduced words}\\
	&\texttt{\indent reduced\_w1=word(p1, commuting\_elements)}\\
	&\texttt{\indent reduced\_w2=word(p2, commuting\_elements)}\\
	&\texttt{\indent \color{red}\#return whether they are equal}\\
	&\texttt{\indent \color{orange}return\color{black}(reduced\_w1==reduced\_w2)}
\end{align*}

\section*{Acknowledgments}
The second author would like to thank the Department of Mathematics at Heriot-Watt University for supporting this summer project. Both authors would like to thank Laura Ciobanu for supervising this project.
\par 
Whilst we have made our best efforts to debug and correctly implement code for this algorithm, there may still be mistakes! If you spot any errors or issues with our code, please let us know.

\newpage
\bibliography{program}
\bibliographystyle{plain}

\uppercase{School of Mathematical and Computer Sciences, Heriot-Watt University, Edinburgh, Scotland, EH14 4AS}
\par 
Email address: \texttt{ggc2000@hw.ac.uk, mj74@hw.ac.uk}

\end{document}